\documentclass{article} % For LaTeX2e
\usepackage{times}
\usepackage{hyperref}
\usepackage{url}
\usepackage{algorithm}
\usepackage{algorithmic}
\usepackage{eqparbox}
 \usepackage{subcaption}
\usepackage{graphicx}
\usepackage{float}

\usepackage{pgfplots} 
\pgfplotsset{compat=newest} 
\pgfplotsset{plot coordinates/math parser=false} 
\newlength\figureheight 
\newlength\figurewidth

\usepackage{amssymb}
\usepackage{amsmath}
\usepackage{amsthm}

\newtheorem{lemma}{Lemma}

\newtheorem{theorem}{Theorem}

\renewcommand{\eqref}[1]{Equation~(\ref{#1})}

\newcommand{\figref}[1]{Figure~\ref{#1}}
\newcommand{\secref}[1]{Section~\ref{#1}}

\newcommand{\appref}[1]{Appendix~\ref{#1}}
\newcommand{\thmref}[1]{Theorem~\ref{#1}}
\newcommand{\lemref}[1]{Lemma~\ref{#1}}

\newcommand{\algref}[1]{Algorithm~\ref{#1}}

\newcommand{\cA}{\mathcal{A}}

\newcommand{\cN}{\mathcal{N}}

\newcommand{\bE}{\mathbb{E}}
\newcommand{\reals}{\mathbb{R}}

\newcommand{\inner}[1]{\langle #1 \rangle}
\newcommand{\tr}{\mathrm{tr}}

\newcommand{\fnorm}[1]{\| #1\|_F}
\newcommand{\snorm}[1]{\| #1\|_{\mathrm{sp}}}
\newcommand{\tnorm}[1]{\| #1\|_{\mathrm{tr}}}

\DeclareMathOperator*{\argmin}{arg\,min}

\newcommand{\w}[1][]{\ifthenelse{\equal{#1}{}}{w}{{w}^{(#1)}} }
\newcommand{\W}[1][]{\ifthenelse{\equal{#1}{}}{W}{{W}^{(#1)}} }
\newcommand{\A}[1][]{\ifthenelse{\equal{#1}{}}{A}{{A}^{(#1)}} }
\newcommand{\vW}[1][]{\ifthenelse{\equal{#1}{}}{W}{{\vec{W}}^{(#1)}} }
\newcommand{\x}[1][]{\ifthenelse{\equal{#1}{}}{x}{{x}^{(#1)}} }
\newcommand{\y}[1][]{\ifthenelse{\equal{#1}{}}{y}{{y}^{(#1)}} }
\newcommand{\g}[1][]{\ifthenelse{\equal{#1}{}}{g}{{g}^{(#1)}} }
\renewcommand{\l}[1][]{\ifthenelse{\equal{#1}{}}{\ell}{{\ell}^{(#1)}} }

\author{Alon Gonen\footnote{School of Computer Science, The Hebrew University, Jerusalem, Isreal} \and Shai Shalev-Shwartz\footnote{School of Computer Science, The Hebrew University, Jerusalem, Isreal}}

\title{Faster SGD Using Sketched Conditioning}

\begin{document}

\maketitle

\begin{abstract}
We propose a novel method for speeding up stochastic optimization algorithms via \emph{sketching} methods, which recently became a powerful tool for accelerating algorithms for numerical linear algebra. We revisit the method of \emph{conditioning} for accelerating first-order methods and suggest the use of sketching methods for constructing a cheap conditioner that attains a significant speedup with respect to the Stochastic Gradient Descent (SGD) algorithm. While our theoretical guarantees assume convexity, we discuss the applicability of our method to deep neural networks, and experimentally demonstrate its merits.
\end{abstract}

\section{Introduction} \label{sec:intro}
We consider empirical loss minimization problems of the form:
\begin{equation} \label{eqn:optProblem}
\min_{W \in \reals^{p \times n}}  L(W) := \frac{1}{m} \sum_{i=1}^m \ell_{y_i}( W x_i ) ~,
\end{equation}
where for every $i$, $\x_i$ is an $n$-dimensional vector and $\ell_{y_i}$ is a loss function from $\reals^p$ to $\reals$.
For example, in multiclass categorization with the logistic-loss, we have that $\ell_{y_i}(a) = \log\left( \sum_{j=1}^p \exp(a_j - a_{y_i})\right)$. 
Later in the paper we will generalize the discussion to the case in which $W$ is the weight matrix of an intermediate layer of a deep neural network. 

We consider the large data regime, in which $m$ is large. A popular algorithm for this case is Stochastic Gradient Descent (SGD). The basic idea is to initialize $W_1$ to be some matrix, and at each time $t$ to draw an index $i \in \{1,\ldots,m\}$ uniformly at random from the training sequence $S=((x_1,y_1),\ldots,(x_m,y_m))$, and then update $W_{t+1}$ based on the gradient of $\ell_{y_i}(W x_i)$ at $W$. When performing this update we would like to decrease the value of $\ell_{y_i}(W x_i)$ while bearing in mind that we only look on a single example, and therefore we should not change $W$ too much. This can be formalized by an update of the form
\[
W_{t+1} = \argmin_{W \in \reals^{p \times n}} \frac{1}{2 \eta} D(W,W_t) + \ell_{y_i}(W x_i) ~,
\]
where $D(\cdot,\cdot)$ is some distance measure between matrices and $\eta$, the \emph{learning rate}, controls the tradeoff between the desire to minimize the function and the desire to stay close to $W_t$. Since we keep $W$ close to $W_t$, we can further simplify things by using the first-order approximation of $\ell_{y_i}$ around $W_t x_i$, 
\[
\ell_{y_i}(W x_i) ~\approx~ \ell_{y_i}(W_t x_i)  + \inner{W -W_t~, ~ (\nabla \ell_{y_i}(W_t x_i) ) x_i^\top } ~,
\]
where $\nabla \ell_{y_i}(W_t x_i) \in \reals^p$ is the (sub)gradient of $\ell_{y_i}$ at the $p$-dimensional vector $W_t x_i$ (as a column vector), and for two matrices $A,B$ we use the notation $\inner{A,B} = \sum_{i,j} A_{i,j} B_{i,j}$. Hence, the update becomes
\begin{equation} \label{eqn:genericPA}
W_{t+1} = \argmin_{W \in \reals^{p \times n}} \frac{1}{2 \eta} D(W,W_t) + \ell_{y_i}(W_t x_i)  + \inner{W -W_t~, ~ (\nabla \ell_{y_i}(W_t x_i) ) x_i^\top }
\end{equation}

\eqref{eqn:genericPA} defines a family of algorithms, where different instances are derived by specifying the distance measure $D$. 
The simplest choice of $D$ is the squared Frobenius norm regularizer, namely, 
\[
D(W,W_t) =  \fnorm{W-W_t}^2 =  \inner{W-W_t,W-W_t} ~.
\]
It is easy to verify that for this choice of $D$, the update given in \eqref{eqn:genericPA} becomes
\[
W_{t+1} = W_t - \eta (\nabla \ell_{y_i}(W_t x_i) ) x_i^\top ~,
\]
which is exactly the update rule of SGD. Note that the Frobenius norm distance measure can be rewritten as
\[
D(W,W_t) = \inner{ I ~,~ (W-W_t)^\top  (W-W_t) }
\]
In this paper, we consider the family of distance measures of the form
\[
D_A(W,W_t) =   \inner{A ~,~ (W-W_t)^\top  (W-W_t) }
\]
where $A$ is a positive definite matrix. For every such choice of $A$, the update given in \eqref{eqn:genericPA} becomes
\begin{equation} \label{eqn:updateGeneric}
W_{t+1} = W_t - \eta (\nabla \ell_{y_i}(W_t x_i) ) (A^{-1} x_i)^\top ~.
\end{equation}
We refer to the matrix $A$ as a conditioning matrix (for a reason that will become clear shortly) and call the resulting algorithm Conditioned SGD.

How should we choose the conditioning matrix $A$? There are two considerations. First, we would like to choose $A$ so that the algorithm will converge to a solution of \eqref{eqn:optProblem} as fast as possible. Second, we would like that it will be easy to compute both the matrix $A$ and the update rule given in \eqref{eqn:updateGeneric}.

We start with the first consideration. Naturally, the convergence of the Conditioned SGD algorithm depends on the specific problem at hand. However, we can rely on convergence bounds and picks $A$ that minimizes these bounds. Concretely, assuming that each $\ell_{y_i}$ is convex and $\rho$-Lipschitz, denote $C = \frac{1}{m} \sum_{i=1}^m x_i x_i^\top$ the correlation matrix of the data, and let $W^*$ be an optimum of \eqref{eqn:optProblem}, then the sub-optimality of the Conditioned SGD algorithm after performing $T$ iterations is upper bounded by 
\[
\frac{1}{2 \eta T} D_A(W^\star, W_1) + \frac{\eta \rho^2}{2} \tr(A^{-1} C) ~.
\]
We still cannot minimize this bound w.r.t. $A$ as we do not know the value of $W^\star$. So, we further upper bound $D_A(W^\star,W_1)$ by considering two possible bounds. Denoting the spectral norm and the trace norm by $\snorm{\cdot}$ and $\tnorm{\cdot}$, respectively, we have
\begin{align*}
&1.~~ D_A(W^\star, W_1) ~\le~ \snorm{A} \, \tnorm{(W^\star-W_1)^\top (W^\star - W_1)} \\
&2.~~ D_A(W^\star, W_1) ~\le~ \tnorm{A} \, \snorm{(W^\star-W_1)^\top (W^\star - W_1)} 
\end{align*}
Interestingly, for the first possibility above, the optimal $A$ becomes $A = I$, corresponding to the vanilla SGD algorithm. However, for the second possibility, we show that the optimal $A$ becomes $A = C^{1/2}$. The ratio between the required number of iterations to achieve $\epsilon$ sub-optimality is
\[
\frac{\textrm{\# iterations for }A = I}{\textrm{\# iterations for }A = C^{1/2}} ~=~ 
\frac{\tnorm{(W^\star-W_1)^\top (W^\star - W_1)} ~~~ \tnorm{C}}{\snorm{(W^\star-W_1)^\top (W^\star - W_1)} ~ \tnorm{C^{1/2}}^2}
\]
The above ratio is always between $1/n$ and $\min\{n,p\}$. We argue that in many typical cases the ratio will be greater than $1$, meaning that the conditioner $A = C^{1/2}$ will lead to faster convergence. For example, suppose that the norm of each row of $W^\star$ is order of $1$, but the rows are not correlated. Let us also choose $W_1 = 0$ and assume that $p = \Theta(n)$. Then, 
$\frac{\tnorm{(W^\star-W_1)^\top (W^\star -
    W_1)}}{\snorm{(W^\star-W_1)^\top (W^\star - W_1)}}$ is order of
$n$. On the other hand, if the eigenvalues of $C$ decay fast, then
$\frac{\tnorm{C}}{\tnorm{C^{1/2}}^2} \approx 1$. Therefore, in such
scenarios, using the conditioner $A = C^{1/2}$ will lead to a factor of $n$ less iterations relatively to vanilla SGD.

Getting back to the question of how to choose $A$, the second consideration that we have mentioned is the time required to compute $A^{-1}$ and to apply the update rule given in \eqref{eqn:updateGeneric}. As we will show later, the time required to compute $A^{-1}$ is less of an issue relatively to the time of applying the update rule at each iteration, so we focus on the latter. 

Observe that the time required to apply \eqref{eqn:updateGeneric} is order of $(p+n)n$. Therefore, if $p \approx n$ then we have no significant overhead in applying the conditioner relatively to applying vanilla SGD. If $p \ll n$, then the update time is dominated by the time required to compute $A^{-1} x_i$. To decrease this time, we propose to use $A$ of the form $Q B Q^\top + a (I - Q Q^\top)$, where $Q \in \reals^{n \times k}$ has orthonormal columns, $B \in \reals^{k \times k}$ is invertible and $k \ll n$. We use linear sketching techniques (see \cite{woodruff2014sketching}) for constructing this conditioner efficiently, and therefore we refer to the resulting algorithm as Sketched Conditioned SGD (SCSGD). Intuitively, the sketched conditioner is a combination of the two conditioners $A = I$ and $A = C^{1/2}$, where the matrix $Q B Q^\top$ captures the top eigenvalues of $C$ and the matrix $a (I - QQ^\top)$ deals with the smaller eigenvalues of $C$. We show that if the eigenvalues of $C$ decay fast enough then SCSGD enjoys similar speedup to the full conditioner $A = C^{1/2}$. The advantage of using the sketched conditioner is that the time required to apply \eqref{eqn:updateGeneric} becomes $(p+k)n$. Therefore, if $p \ge k$ then the runtime per iteration of SCSGD and the runtime per iteration of the vanilla SGD are of the same order of magnitude.

The rest of the paper is organized as follows. In the next subsection
we survey some related work. In \secref{sec:cond} we describe in detail our conditioning method.
Finally, in \secref{sec:experiments} we discuss variants of the method that are applicable to deep
learning problems and report some preliminary
experiments showing the merits of conditioning for deep learning problems. Due to the lack of space, proofs are
omitted and can be found in \appref{app:proofs}.

\subsection{Related work} \label{sec:related}
Conditioning is a well established technique in optimization aiming at
choosing an ``appropriate'' coordinate system for the optimization process. 
For twice differentiable objectives, maybe the most well known
approach is Newton's method which dynamically changes the coordinate
system according to the Hessian of the objective around the current
solution. There are several problems with utilizing the
Hessian. First, in our case, the Hessian matrix is of size $(pn)
\times (pn)$. Hence, it is computationally expensive to compute and
invert it. Second, even for convex problems, the Hessian matrix might
be meaningless. For example, for linear regression with the absolute
loss the Hessian matrix is the zero matrix almost everywhere. Third,
when the number of training examples is very large, stochastic methods
are preferable and it is not clear how to adapt Newton method to the
stochastic case. The crux of the problem is that while it is easy to
construct an unbiased, low variance, estimate of the gradient, based
on a single example, it is not clear how to construct a good estimate
of the Newton's direction based on a small mini-batch of examples. 

Many approaches have been proposed for speeding up Newton's method.
For example, the $R\{\cdot\}$ operator technique
\cite{pearlmutter1994fast,werbos1988backpropagation,moller1993exact,martens2010deep}. However,
these methods are not applicable for the stochastic setting.  An
obvious way to decrease the storage and computational cost is to only
consider the diagonal elements of the Hessian (see
\cite{becker1988improving}).  Schraudolph
\cite{schraudolph2007stochastic} proposed an adaptation of the L-BFGS
approach to the online setting, in which at each iteration, the
estimation of the inverse of the Hessian is computed based on only the
last few noisy gradients. Naturally, this yields a low rank
approximation. In \cite{bordes2009sgd}, the two aforementioned
approaches are combined to yield the SGD-QN algorithm.  In the same
paper, an analysis of second order SGD is described, but with $A$
being always the Hessian matrix at the optimum (which is of course not
known). There are various other approximations, see for example
\cite{saito1997partial,bengio2012practical,vinyals2011krylov,roux2008topmoumoute}.

To tackle the second problem, several methods
\cite{schraudolph2002fast,martens2010deep,vinyals2011krylov,pascanu2013revisiting}
rely on different variants of the Gauss-Newton approximation of the
Hessian. A somewhat related approach is Amari's natural gradient
descent \cite{amari1998natural,amari2000adaptive}. See the discussion
in \cite{pascanu2013revisiting}. To the best of our knowledge, these
methods come with no theoretical guarantees. 

The aforementioned approaches change the conditioner at each iteration
of the algorithm. A general treatment of this approach is described in 
\cite{nesterov2004introductory}[Section 1.3.1] under the name
``Variable Metric''. 
Maybe the most relevant approach is the Adagrad algorithm
\cite{duchi2011adaptive}, which was originally proposed for the online
learning setting but can be easily adapted to the stochastic
optimization setting. In our notation, the AdaGrad algorithm uses 
a $(pn) \times (pn)$ conditioning matrix that changes along time and has the form, 
$A_t = \delta I + \frac{1}{t} \sum_{i=1}^t
\nabla_t \nabla_t^\top $, where $\nabla_t = \mathrm{vec}(\nabla
\ell_{y_i}(W_t x_i) ) x_i^\top)$.
There are several advantages of our method relatively to
AdaGrad. First, the convergence bound we obtain is better than the
convergence bound of AdaGrad. Specifically, while both bounds have the
sane dependence on $\tnorm{C^{1/2}}^2$, our bound depends on
$\snorm{W^*}^2$ while AdaGrad depends on $\fnorm{W^*}^2$. As we discussed
before, there may be a gap of $p$ between $\snorm{W^*}^2$ and
$\fnorm{W^*}^2$. More critically, when using a full conditioner, the
storage complexity of our conditioner is $n^2$, while the storage
complexity of AdaGrad is $(np)^2$. In addition, the time complexity of
applying the update rule is $(p+n)n$ for our conditioner versus
$(np)^2$ for AdaGrad. For this reason, most practical application of
AdaGrad relies on a diagonal approximation of $A_t$. In contrast, we
can use a full conditioner in many practical cases, and even when $n$
is large our sketched conditioner can be applied without a significant
increase in the complexity relatively to vanilla SGD. 
Finally, because we derive our algorithm for the stochastic case (as
opposed to the adversarial online optimization setting), and because
we bound the component $\nabla \ell_y (W_t x)$ using the Lipschitzness
of $\ell_y$, the conditioner we use is the constant $C^{-1/2}$ along
the entire run of the optimization process, and should only be
calculated once. In contrast, AdaGrad replaces the conditioner in every
iteration.

\section{Conditioning and Sketched Conditioning} \label{sec:cond}
As mentioned previously, the algorithms we consider start with an
initial matrix $W_1$ and at each iteration update the matrix
according to \eqref{eqn:updateGeneric}. The following lemma provides
an upper bound on the expected sub-optimality of any algorithm of this
form. 
 \begin{lemma} \label{lem:md}
Fix a positive definite matrix $A \in \reals^{n \times n}$. Let
$W^\star$ be the minimizer of \eqref{eqn:optProblem}, let $\sigma \in
\reals$ be such that $\sigma \ge \snorm{W^\star}$ and denote $C =
\frac{1}{m}\sum_{i=1}^m x_i x_i^\top$. Assume that for every $i$,
$\ell_{y_i}$ is convex and $\rho$-Lipschitz. Then, if we apply the
update rule given in \eqref{eqn:updateGeneric} using the conditioner $A$ and denote $\bar{W} =
\frac{1}{T} \sum_{t=1}^T W_t$, then 
\begin{align*}
\bE[L(\bar{W})- L(W^\star) ] &\le \frac{1}{2\eta T}  \tr(A {W^\star}^\top W^\star)  + \frac{\eta \rho^2}{2} \bE \left[ \tr(A^{-1} C) \right] \\&
\le \frac{\sigma^2}{2\eta T} \tr(A)  + \frac{\eta \rho^2}{2} \bE \left[ \tr(A^{-1} C) \right]~.
\end{align*}
In particular, for $\eta= \sigma/(\rho \sqrt{T})$, we obtain
\[
\bE[L(\bar{W})- L(W^\star) ] \le  \frac{\sigma \rho}{\sqrt{T}} (\tr(A) + \tr(A^{-1}C))~.
\]
\end{lemma}
The proof of the above lemma can be obtained by replacing the standard inner product with the inner product induced by $A$. For
completeness, we provide a proof in \appref{app:proofs}.

In \appref{app:proofs} we show that the conditioner which minimizes
the bound given in the above Lemma is $A = C^{1/2}$. This yields:
\begin{theorem} \label{thm:optBound}
Following the notation of \lemref{lem:md}, assume that we run the
meta-algorithm with $A = C^{1/2}$, then
\[
 \bE[L(\bar{W}) - L(W^\star)]   \le \frac{\sigma \rho}{\sqrt{T}} \cdot \tr(C^{1/2})~.
\]
\end{theorem}

\subsection{Sketched Conditioning}
Let $k < n$ and assume that $\textrm{rank}(C) \ge k$. Consider the following family of conditioners:
\begin{equation} \label{eq:lowCond}
{\cA} = \{A = Q B Q^\top + a (I - Q Q^\top): Q \in \reals^{n \times k},\, Q^\top Q=I,\,B \succ 0 \in \reals^{k \times k},\, a>0 \}
\end{equation}
Before proceeding, we show that the conditioners in $\cA$ are indeed
positive definite, and give a formula for their inverse.
\begin{lemma} \label{lem:lowCondInv}
Let $A = Q B Q^\top + a (I-QQ^\top) \in A$. Then, $A \succ 0$ and its inverse is given by
\[
{A}^{-1} = Q B^{-1}  Q^\top + a^{-1} (I-QQ^\top)~.
\]
\end{lemma}
Informally, every conditioner $A \in \cA$ is a combination of a low
rank conditioner and the identity conditioner. The most appealing
property of these conditioners is that we can compute 
$A^{-1} x$ in time $O(nk)$ and therefore the time complexity of
calculating the update given in \eqref{eqn:updateGeneric} is $O(n(p+k))$.

In the next subsections we focus on instances of $\cA$ which are induced by an approximate best rank-$k$ approximation of $C$. However, for now, we give an analysis for any choice of $A \in \cA$.  
\begin{theorem} \label{thm:generalLowCond}
Following the notation of \lemref{lem:md}, let $A \in \cA$ and denote $\tilde{C}=Q^\top C Q$. Then, if $a = \sqrt{\frac{\tr(C)-\tr(\tilde{C})}{n-k}}$, we have
\[
\bE[L(\bar{W})- L(W^\star) ] \le \frac{\sigma \rho}{2\sqrt{T}} \cdot \left(\tr(B) + \tr(B^{-1} \tilde{C}) + 2\sqrt{(n-k)(\tr(C) - \tr(\tilde{C}))} \right)~.
\]
\end{theorem}

\subsection{Low-rank conditioning via exact low-rank
  approximation} \label{sec:condExact} Maybe the most straightforward
approach of defining $Q$ and $B$ is by taking the leading eigenvectors
of $C$. Concretely, let $C=UDU^\top$ be the eigenvalue decomposition
of $C$ and denote the diagonal elements of $D$ by $\lambda_1 \ge
\ldots \ge \lambda_n \ge 0$. Recall that for any $k \le n$, the best
rank-$k$ approximation of $C$ is given by $C_k = U_k D_k U_k^\top$,
where $U_k \in \reals^{n \times k}$ consists of the first $k$ columns of $U$
and $D_k$ is the first $k \times k$ sub-matrix of $D$.  Denote
$\tilde{C} = Q^\top C Q$ and consider the conditioner $\tilde{A}$
which is determined from \eqref{eq:lowCond} by setting $Q = U_k$, $B =
\tilde{C}^{1/2}$, and $a$ as in \thmref{thm:generalLowCond}. 

\begin{theorem} \label{thm:lowCondExact}
Let $Q = U_k$, $B =
\tilde{C}^{1/2}$ and $a$ as in \thmref{thm:generalLowCond},  and
consider the conditioner given in \eqref{eq:lowCond}. Then, 
\[
\bE[L(\bar{W})- L(W^\star) ] \le \frac{\sigma \rho}{\sqrt{T}} \cdot \left(\tr(C_k^{1/2}) + \sqrt{(n-k) (\tr(C)-\tr(C_{k}))} \right)~.
\]
In particular, if $\sqrt{(n-k)  (\tr(C)-\tr(C_{k}))} =
O(\tr(C^{1/2}))$, then the obtained bound is of the same order as the bound in \thmref{thm:optBound}.
\end{theorem}
We refer to the condition $\sqrt{(n-k)  (\tr(C)-\tr(C_{k}))} =
O(\tr(C^{1/2}))$ as a \emph{fast spectral decay property} of the matrix $C$.

\subsection{Low-rank conditioning via sketching}  \label{sec:condSketch}
The conditioner defined in the previous subsection requires
the exact calculation of the matrix $C$ and its eigenvalue
decomposition. In this section we describe a faster technique for
calculating a sketched conditioner. Before formally describing the sketching technique, let us try to explain the intuition behind it. \figref{fig:jl} depicts a set of 1000 (blue) random points in the plane. Suppose that we represent this sequence by a matrix $X \in \reals^{2 \times 1000}$. Now we draw a vector $\omega \in \reals^{1000 \times 1}$ whose coordinates are $\cN(0,1)$ i.i.d. random variables and consider the vector $z=X\omega$. The vector $z$ is simply a random combination of these points. As we can see, $z$ coincides with the strongest direction of the data. More generally, the idea of sketching is that if we take a matrix $X \in \reals^{n \times m}$ and multiply it from the right by random matrix $\Omega \in \reals^{m \times r}$, then with high probability, we preserve the strongest directions of the column space of $X$. The above intuition is formalized by the following result, which follows from \cite{sarlos2006improved} by setting $\epsilon=1$\footnote{See also
  \cite{woodruff2014sketching}[Lemmas 4.1,4.2]. In particular, the
  elements of $\Omega$ can be drawn either according to be
  i.i.d. $\cN(0,1/r)$ or zero-mean $\pm 1$ random variables. Also, the
  bounds on the lower dimension in \cite{woodruff2014sketching} are
  better in (additive) factor $k \log k$.}.

\begin{figure}
\centering
\input{JL.tex}
\caption{}
\label{fig:jl}
\end{figure}

\begin{lemma}  \label{lem:woodruff}
Let $X \in \reals^{n \times m}$. Let $r = \Theta(k)$ and let $\Omega \in \reals^{m \times r}$ be a random matrix whose elements are i.i.d. $\cN(0,1/r)$ random variables. 
Let $P \in \reals^{n \times r}$ be a matrix whose columns form an orthonormal basis of the column space of $X \Omega$, 
let $U \in \reals^{r \times k}$ be a matrix whose columns are the top $k$ eigenvectors of the matrix $(P^\top X)(P^\top X)^\top$, and 
let $Q = P U \in \reals^{n \times k}$. Then, 
\begin{align} \label{eq:woodruff}
\bE\|QQ^\top X - X\|_F \le 2 \|X-X_k\|_F~.
\end{align}
\end{lemma}
Let $X \in \reals^{ n \times m}$ be a matrix whose columns are the
vectors $x_1,\ldots,x_m$. Based on \lemref{lem:woodruff}, we produce a
matrix $Q \in \reals^{n \times k}$ which satisfies the inequality
$\bE[\|QQ^\top X - X\|_F ]\le 2 \|X-X_k\|_F$. Let $\tilde{C} = Q^\top C Q$. Our sketched conditioner is determined by the matrix $Q$ and the matrix $B=\tilde{C}^{1/2}$. As we show in \algref{alg:condSketch}, we can compute a factored form of the inverse of the conditioner, $\tilde{A}^{-1}$, in time $O(mnk)$. 
\begin{algorithm}
\caption{Sketched Conditioning: Preprocessing}
\label{alg:condSketch}
\begin{algorithmic}
\STATE \textbf{Input: } $X \in \reals^{n,m}$ ~~~,~~~ \textbf{Parameters: } $k < n, r \in \Theta(k)$ 
\STATE \textbf{Output: }
$Q$, $B^{-1}$, $a^{-1}$ that determines a conditioner according to \eqref{eq:lowCond}
\STATE Sample each element of $\Omega \in \reals^{m \times r}$ i.i.d. from
$\cN(0,r^{-1})$
\STATE Compute $Z = X \Omega$ \COMMENT {in time  $O(mnr)$}
\STATE $[P,\sim] = \mathrm{QR}(Z)$   \COMMENT {in time $O(r^2n)$}
\STATE Compute $Y = P^\top X$ \COMMENT {in time $O(mnr)$}
\STATE Compute the SVD: $~~ Y = U' \Sigma' V'^\top$ \COMMENT {in time $O(mr^2)$}
\STATE Compute $Q = PU'_k$ \COMMENT{in time $O(nrk)$}
\STATE Compute $\tilde{C}=Q^\top C Q$  \COMMENT{in time $O(mkn)$} 
\STATE Compute $B^{-1} = \tilde{C}^{-1/2}$ \COMMENT{in time $O(k^3)$} 
\STATE Compute $a^{-1} = \sqrt{\frac{n-k} {\tr(C)-\tr(\tilde{C})}}$  \COMMENT{in time $O(mn+k)$}
\end{algorithmic}
\end{algorithm}
We turn to discuss the performance of this conditioner. Relying on \lemref{lem:woodruff}, we start by relating the trace of $\tilde{C}=Q^\top C Q$ to the trace of $C$.
\begin{lemma}  \label{lem:tildeC_C}
We have $\tr(C) - \tr(\tilde{C}) \le 4 (\tr(C)-\tr(C_k))$.
\end{lemma}
The next lemma holds for any choice of $Q \in \reals^{n \times k}$ with orthonormal columns. 
\begin{lemma} \label{lem:everyQ}
Assume that $C$ is of rank at least $k$. Let $Q \in \reals^{n \times k}$ with orthonormal columns and define $\tilde{C} = Q^\top C Q^\top$, $B=\tilde{C}^{1/2}$. Then, $\tr(B) = \tr(B^{-1} \tilde{C})= O( \tr(C_k^{1/2}))$.

\end{lemma}
Combining the last two lemmas with \thmref{thm:generalLowCond}, we conclude:
\begin{theorem} \label{cor:alonHofer}
Consider running SCSGD with the conditioner given in \algref{alg:condSketch}. Then,
\begin{align*}
\bE[L(\bar{W})- L(W^\star) ] &\le O\left( \frac{\sigma \rho}{\sqrt{T}} \cdot \left(\tr(C_k^{1/2}) + \sqrt{(n-k) (\tr(C)-\tr(C_k))} \right) \right)~.
\end{align*}
In particular, if the fast spectral decay property holds, i.e., $\sqrt{(n-k)(\tr(C)-\tr(C_k))} = O(\tr(C^{1/2}))$, then the obtained bound is of the same order as the bound in \thmref{thm:optBound}.
\end{theorem}

\section{Experiments with Deep Learning} \label{sec:experiments} 
While our theoretical guarantees were derived for convex problems, the
conditioning technique can be adapted for deep learning problems, as
we outline below.

A feedforward deep neural network is a function $f$ that can be
written as a composition $f = f_1 \circ f_2 \circ \ldots \circ f_q$,
where each $f_i$ is called a layer function. Some of the layer
functions are predefined, while other layer functions are
parameterized by weights matrices. Training of a network amounts
to optimizing w.r.t. the weights matrices. The most popular layer
function with weights is the affine layer (a.k.a. ``fully connected''
layer). This layer performs the transformation $y = W x + b$, where $x
\in \reals^n$, $W \in \reals^{p,n}$, and $b \in \reals^p$. The network
is usually trained based on variants of stochastic gradient descent,
where the gradient of the objective w.r.t. $W$ is calculated based on
the backpropagation algorithm, and has the form $\delta x^\top$, where
$\delta \in \reals^p$.

To apply conditioning to an affine layer, 
instead of the vanilla SGD update $W = W - \eta \delta x^\top$, we can
apply a conditioned update of the form $W = W - \eta \delta (A^{-1 }
x)^\top$. To calculate $A$ we could go over the entire training data
and calculate $C = \frac{1}{m} \sum_{i=1}^m x_i x_i^\top$. However,
unlike the convex case, now the vectors $x_i$ are not constant but
depends on weights of previous layers. Therefore, we initialize $C =
I$ and update it according to the update rule $C = (1-\nu) C + \nu x_i x_i^\top$.
for some $\nu \in (0,1)$. From time to time, we replace the
conditioner to be $A = C^{1/2}$ for the current value of $A$. In our
experiments, we updated the conditioning matrix after each $50$s
iterations. Note that the process of calculating $A = C^{1/2}$ can be
performed in a different thread, in parallel to the main stochastic
gradient descent process, and therefore it causes no slowdown to the main
stochastic gradient descent process.

The same technique can be applied to convolutional layers (that also
have weights), because it is possible to write a convolutional layer
as a composition of a transformation called ``Im2Col'' and a vanilla
affine layer. Besides these changes, the rest of the algorithm is the same as in the
convex case.
 
Below we describe two experiments in which we have applied
conditioning technique to a popular variant of stochastic gradient
descent. In particular, we used stochastic gradient descent with a
mini-batch of size $64$, a learning rate of $\eta_t = 0.01 (1 + 0.0001
t)^{-3/4}$, and with Nesterov momentum with parameter $0.9$, as
described in \cite{sutskever2013importance}. To initialize the weights
we used the so-called Xavier method, namely, chose each element of $W$
at random according to a uniform distribution over $[-a,a]$, with $a =
\sqrt{3/n}$. We chose these parameters because they are the default in
the popular Caffe library (\url{http://caffe.berkeleyvision.org}),
without attempting to tune them. We conducted experiments with the
MNIST dataset \cite{lecun1998mnist} and with the Street View House
Numbers (SVHN) dataset \cite{netzer2011reading}.

\paragraph{MNIST: } We used a variant of the LeNet architecture
\cite{lecun1998gradient}. The input is images of $28 \times 28$ pixels. We
apply the following layer functions: Convolution with kernel size of
$5 \times 5$, without padding, and with 20 output
channels. Max-pooling with kernel size of $2 \times 2$. Again, a
convolutional and pooling layers with the same kernel sizes, this time
with 50 output channels. Finally, an affine layer with 500 output
channels, followed by a ReLU layer and another affine layer with 10
output channels that forms the prediction. In short, the architecture
is: conv 5x5x20, maxpool 2x2, conv 5x5x50, maxpool 2x2, affine 500,
relu, affine 10.

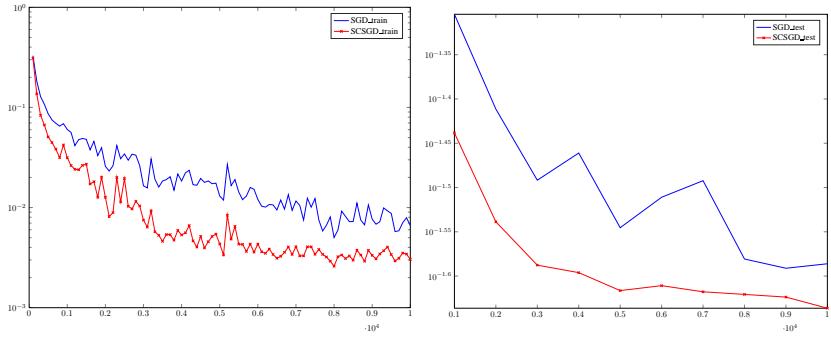
\begin{figure}
\begin{subfigure} [b] {0.45 \textwidth}
\centering
\resizebox{\linewidth}{!}{% This file was created by matlab2tikz.
% Minimal pgfplots version: 1.3
%  http://www.mathworks.com/matlabcentral/fileexchange/22022-matlab2tikz
\begin{tikzpicture}

\begin{axis}[%
width=6.027778in,
height=4.754167in,
at={(1.011111in,0.641667in)},
scale only axis,
xmin=0,
xmax=10000,
ymode=log,
ymin=0.001,
ymax=1,
yminorticks=true,
legend style={legend cell align=left,align=left,draw=white!15!black}
]
\addplot [color=blue,solid]
  table[row sep=crcr]{%
100	0.328427\\
200	0.183425\\
300	0.128119\\
400	0.107447\\
500	0.0862663\\
600	0.0750614\\
700	0.069434\\
800	0.0650359\\
900	0.0688158\\
1000	0.0602614\\
1100	0.0562707\\
1200	0.0415663\\
1300	0.0476752\\
1400	0.0489894\\
1500	0.0479795\\
1600	0.0377317\\
1700	0.0458077\\
1800	0.0330539\\
1900	0.0396087\\
2000	0.0258173\\
2100	0.0231698\\
2200	0.0261381\\
2300	0.0415923\\
2400	0.030713\\
2500	0.0343127\\
2600	0.0296985\\
2700	0.0341199\\
2800	0.03334\\
2900	0.0263124\\
3000	0.0164864\\
3100	0.015723\\
3200	0.0304979\\
3300	0.0191931\\
3400	0.0160307\\
3500	0.0184158\\
3600	0.0190374\\
3700	0.0203258\\
3800	0.0149335\\
3900	0.0216347\\
4000	0.0184254\\
4100	0.0222262\\
4200	0.0235551\\
4300	0.0169395\\
4400	0.0167014\\
4500	0.0195337\\
4600	0.0178428\\
4700	0.0184471\\
4800	0.0173052\\
4900	0.017542\\
5000	0.0130179\\
5100	0.0118754\\
5200	0.0269306\\
5300	0.0165976\\
5400	0.0190315\\
5500	0.0142642\\
5600	0.0119987\\
5700	0.0130323\\
5800	0.0158637\\
5900	0.0152257\\
6000	0.012062\\
6100	0.0103223\\
6200	0.0101095\\
6300	0.0107368\\
6400	0.0106317\\
6500	0.00946309\\
6600	0.0118979\\
6700	0.00965037\\
6800	0.0134538\\
6900	0.00939099\\
7000	0.0116149\\
7100	0.0104989\\
7200	0.00749987\\
7300	0.0123951\\
7400	0.0100903\\
7500	0.0122955\\
7600	0.00751995\\
7700	0.00583415\\
7800	0.00666049\\
7900	0.00806003\\
8000	0.00502859\\
8100	0.00594521\\
8200	0.00920762\\
8300	0.00815611\\
8400	0.00724418\\
8500	0.00724964\\
8600	0.0110521\\
8700	0.0075117\\
8800	0.00671469\\
8900	0.0106343\\
9000	0.00766212\\
9100	0.00682251\\
9200	0.00721176\\
9300	0.00992564\\
9400	0.00927391\\
9500	0.00874279\\
9600	0.00575291\\
9700	0.0058661\\
9800	0.0070488\\
9900	0.00791083\\
10000	0.00660616\\
};
\addlegendentry{SGD\_train};

\addplot [color=red,solid,mark=x,mark options={solid}]
  table[row sep=crcr]{%
100	0.309666\\
200	0.136498\\
300	0.0833322\\
400	0.0668218\\
500	0.050764\\
600	0.0445666\\
700	0.0383272\\
800	0.0315715\\
900	0.0421541\\
1000	0.0314995\\
1100	0.0262857\\
1200	0.0241018\\
1300	0.0238968\\
1400	0.0263785\\
1500	0.0271212\\
1600	0.0172008\\
1700	0.018147\\
1800	0.012703\\
1900	0.0200644\\
2000	0.0126821\\
2100	0.00812317\\
2200	0.00891333\\
2300	0.0199609\\
2400	0.0114035\\
2500	0.0195042\\
2600	0.0103019\\
2700	0.009694\\
2800	0.0115497\\
2900	0.0103847\\
3000	0.00748938\\
3100	0.00638315\\
3200	0.00928139\\
3300	0.00573185\\
3400	0.00528243\\
3500	0.00460685\\
3600	0.00537872\\
3700	0.00536445\\
3800	0.00472765\\
3900	0.00591746\\
4000	0.00532436\\
4100	0.00559794\\
4200	0.00660649\\
4300	0.00465408\\
4400	0.00403844\\
4500	0.00514708\\
4600	0.00395344\\
4700	0.00456138\\
4800	0.00514766\\
4900	0.00543171\\
5000	0.00432912\\
5100	0.00337217\\
5200	0.00838169\\
5300	0.00485775\\
5400	0.00647529\\
5500	0.00431119\\
5600	0.00428784\\
5700	0.00364759\\
5800	0.00431998\\
5900	0.00359191\\
6000	0.0043182\\
6100	0.00362031\\
6200	0.00349343\\
6300	0.00384749\\
6400	0.00340223\\
6500	0.00311979\\
6600	0.00326216\\
6700	0.00357251\\
6800	0.00404201\\
6900	0.00340662\\
7000	0.00406575\\
7100	0.0032917\\
7200	0.00330025\\
7300	0.00405448\\
7400	0.00405619\\
7500	0.00341556\\
7600	0.003817\\
7700	0.00340697\\
7800	0.00320825\\
7900	0.00291741\\
8000	0.00259668\\
8100	0.00322956\\
8200	0.00335876\\
8300	0.00309913\\
8400	0.00328742\\
8500	0.00298466\\
8600	0.00375257\\
8700	0.00335669\\
8800	0.00291301\\
8900	0.00374011\\
9000	0.00333371\\
9100	0.00307392\\
9200	0.00344784\\
9300	0.00370976\\
9400	0.00402974\\
9500	0.00339907\\
9600	0.00292968\\
9700	0.00311563\\
9800	0.00351304\\
9900	0.00343025\\
10000	0.00304097\\
};
\addlegendentry{SCSGD\_train};

\end{axis}
\end{tikzpicture}%
}
\end{subfigure}
\begin{subfigure} [b] {0.45 \textwidth}
\centering
\resizebox{\linewidth}{!}{% This file was created by matlab2tikz.
% Minimal pgfplots version: 1.3
%  http://www.mathworks.com/matlabcentral/fileexchange/22022-matlab2tikz
\begin{tikzpicture}

\begin{axis}[%
width=6.027778in,
height=4.754167in,
at={(1.011111in,0.641667in)},
scale only axis,
xmin=1000,
xmax=10000,
ymode=log,
ymin=0.0230993,
ymax=0.0496065,
yminorticks=true,
legend style={legend cell align=left,align=left,draw=white!15!black}
]
\addplot [color=blue,solid]
  table[row sep=crcr]{%
1000	0.0496065\\
2000	0.0387792\\
3000	0.0322238\\
4000	0.034582\\
5000	0.0284792\\
6000	0.0308258\\
7000	0.0321838\\
8000	0.0262519\\
9000	0.0256294\\
10000	0.0259303\\
};
\addlegendentry{SGD\_test};

\addplot [color=red,solid,mark=x,mark options={solid}]
  table[row sep=crcr]{%
1000	0.0364411\\
2000	0.0289206\\
3000	0.0258375\\
4000	0.0253388\\
5000	0.0241826\\
6000	0.0244943\\
7000	0.0241071\\
8000	0.0239436\\
9000	0.0237776\\
10000	0.0230993\\
};
\addlegendentry{SCSGD\_test};

\end{axis}
\end{tikzpicture}%
%%% Local Variables:
%%% mode: latex
%%% TeX-master: t
%%% End:
}
\end{subfigure}
\caption{MNIST data. Train (left) and test (right) errors w.r.t. the multiclass log loss of SGD and SCSGD}
\label{fig:mnist}
\end{figure}

\begin{figure}
\begin{subfigure} [b] {0.45 \textwidth}
\centering
\resizebox{\linewidth}{!}{% This file was created by matlab2tikz.
% Minimal pgfplots version: 1.3
%  http://www.mathworks.com/matlabcentral/fileexchange/22022-matlab2tikz
\begin{tikzpicture}

\begin{axis}[%
width=6.027778in,
height=4.754167in,
at={(1.011111in,0.641667in)},
scale only axis,
xmin=0,
xmax=10000,
ymin=0,
ymax=0.1,
legend style={legend cell align=left,align=left,draw=white!15!black}
]
\addplot [color=blue,solid]
  table[row sep=crcr]{%
100	0.0997818\\
200	0.0541456\\
300	0.0396486\\
400	0.0313273\\
500	0.0281567\\
600	0.0208609\\
700	0.0206755\\
800	0.0230302\\
900	0.0207858\\
1000	0.0171644\\
1100	0.0175065\\
1200	0.0111245\\
1300	0.0145806\\
1400	0.0119627\\
1500	0.0121788\\
1600	0.0130676\\
1700	0.0128475\\
1800	0.0111724\\
1900	0.0106431\\
2000	0.00726477\\
2100	0.00792237\\
2200	0.00830534\\
2300	0.0108716\\
2400	0.00760803\\
2500	0.0110031\\
2600	0.00946071\\
2700	0.0101793\\
2800	0.0091882\\
2900	0.0092131\\
3000	0.00483167\\
3100	0.00354433\\
3200	0.0078924\\
3300	0.00660888\\
3400	0.0037256\\
3500	0.00478963\\
3600	0.00709902\\
3700	0.00842687\\
3800	0.00436074\\
3900	0.0079564\\
4000	0.00624883\\
4100	0.00644663\\
4200	0.00690099\\
4300	0.00498\\
4400	0.00374072\\
4500	0.00426886\\
4600	0.00635266\\
4700	0.00648746\\
4800	0.00538409\\
4900	0.0051769\\
5000	0.00262828\\
5100	0.00306523\\
5200	0.00745872\\
5300	0.00465254\\
5400	0.00552681\\
5500	0.00325151\\
5600	0.00326043\\
5700	0.00397376\\
5800	0.00510528\\
5900	0.00328073\\
6000	0.00321905\\
6100	0.00339661\\
6200	0.00310651\\
6300	0.00298233\\
6400	0.00350003\\
6500	0.00216939\\
6600	0.00295917\\
6700	0.00217183\\
6800	0.00346487\\
6900	0.00307526\\
7000	0.00287396\\
7100	0.00302668\\
7200	0.00181818\\
7300	0.0022921\\
7400	0.00291846\\
7500	0.00380007\\
7600	0.00183436\\
7700	0.000522856\\
7800	0.0010889\\
7900	0.00141809\\
8000	0.00123165\\
8100	0.001274\\
8200	0.00159254\\
8300	0.00261606\\
8400	0.0020387\\
8500	0.00182638\\
8600	0.00247101\\
8700	0.00175445\\
8800	0.00148715\\
8900	0.00270461\\
9000	0.00227381\\
9100	0.00087093\\
9200	0.00235856\\
9300	0.0025758\\
9400	0.00275676\\
9500	0.0025316\\
9600	0.000473217\\
9700	0.000942792\\
9800	0.00101763\\
9900	0.0017767\\
10000	0.000718144\\
};
\addlegendentry{\large{SGD\_train\_01}};

\addplot [color=red,solid,mark=x,mark options={solid}]
  table[row sep=crcr]{%
100	0.0960575\\
200	0.0404565\\
300	0.0222333\\
400	0.0177711\\
500	0.013487\\
600	0.0122223\\
700	0.00983849\\
800	0.011707\\
900	0.0117181\\
1000	0.0087044\\
1100	0.00718554\\
1200	0.00714726\\
1300	0.00729439\\
1400	0.00642668\\
1500	0.0048555\\
1600	0.00345858\\
1700	0.00574181\\
1800	0.0029268\\
1900	0.00527772\\
2000	0.00281503\\
2100	0.00135752\\
2200	0.00163639\\
2300	0.00402942\\
2400	0.00164102\\
2500	0.00378628\\
2600	0.00143156\\
2700	0.0021195\\
2800	0.00247575\\
2900	0.00262351\\
3000	0.00203805\\
3100	0.000647972\\
3200	0.00272835\\
3300	0.00050019\\
3400	0.000402551\\
3500	0.000277791\\
3600	0.00013765\\
3700	0.000862153\\
3800	0.000334019\\
3900	0.000274569\\
4000	0.000551197\\
4100	0.000116366\\
4200	0.000421186\\
4300	4.17932e-05\\
4400	0.000385495\\
4500	0.000151166\\
4600	8.94983e-07\\
4700	5.29877e-09\\
4800	0.000299814\\
4900	0.00125611\\
5000	0.000446\\
5100	2.53245e-05\\
5200	0.00158579\\
5300	0.000390341\\
5400	0.000403111\\
5500	7.25548e-06\\
5600	0.00029485\\
5700	7.14052e-06\\
5800	3.18945e-05\\
5900	2.28663e-05\\
6000	1.35381e-07\\
6100	1.16454e-05\\
6200	9.5449e-05\\
6300	0.000781815\\
6400	4.88175e-05\\
6500	2.89025e-07\\
6600	6.32788e-05\\
6700	3.74644e-07\\
6800	0.000762997\\
6900	4.51734e-06\\
7000	0.000753004\\
7100	4.45818e-06\\
7200	2.64838e-05\\
7300	9.55369e-05\\
7400	0.000401614\\
7500	1.13881e-05\\
7600	6.74238e-08\\
7700	3.99184e-10\\
7800	2.36338e-12\\
7900	1.39924e-14\\
8000	8.28424e-17\\
8100	4.90471e-19\\
8200	2.90384e-21\\
8300	1.71923e-23\\
8400	1.01787e-25\\
8500	8.55986e-06\\
8600	5.06788e-08\\
8700	3.00045e-10\\
8800	1.77642e-12\\
8900	2.64574e-05\\
9000	1.56642e-07\\
9100	9.27401e-10\\
9200	5.4907e-12\\
9300	2.81975e-05\\
9400	0.000669991\\
9500	3.96669e-06\\
9600	7.74876e-06\\
9700	4.58767e-08\\
9800	3.59923e-05\\
9900	3.62051e-05\\
10000	2.14353e-07\\
};
\addlegendentry{\large{SCSGD\_train\_01}};

\end{axis}
\end{tikzpicture}%
%%% Local Variables:
%%% mode: latex
%%% TeX-master: t
%%% End:
}
\end{subfigure}
\begin{subfigure} [b] {0.45 \textwidth}
\centering
\resizebox{\linewidth}{!}{% This file was created by matlab2tikz.
% Minimal pgfplots version: 1.3
%  http://www.mathworks.com/matlabcentral/fileexchange/22022-matlab2tikz
\begin{tikzpicture}

\begin{axis}[%
width=6.027778in,
height=4.754167in,
at={(1.011111in,0.641667in)},
scale only axis,
xmin=1000,
xmax=10000,
ymin=0.007,
ymax=0.016,
legend style={legend cell align=left,align=left,draw=white!15!black}
]
\addplot [color=blue,solid]
  table[row sep=crcr]{%
1000	0.0159\\
2000	0.0127\\
3000	0.0102\\
4000	0.0111\\
5000	0.0095\\
6000	0.0103\\
7000	0.0112\\
8000	0.0084\\
9000	0.0087\\
10000	0.0076\\
};
\addlegendentry{\large{SGD\_test\_01}};

\addplot [color=red,solid,mark=x,mark options={solid}]
  table[row sep=crcr]{%
1000	0.0123\\
2000	0.009\\
3000	0.0078\\
4000	0.008\\
5000	0.0077\\
6000	0.0081\\
7000	0.0084\\
8000	0.008\\
9000	0.008\\
10000	0.0075\\
};
\addlegendentry{\large{SCSGD\_test\_01}};

\end{axis}
\end{tikzpicture}%
%%% Local Variables:
%%% mode: latex
%%% TeX-master: t
%%% End:
}
\end{subfigure}
\caption{MNIST data. Train (left) and test (right) errors w.r.t. the zero-one loss of SGD and SCSGD}
\label{fig:mnist01}
\end{figure}
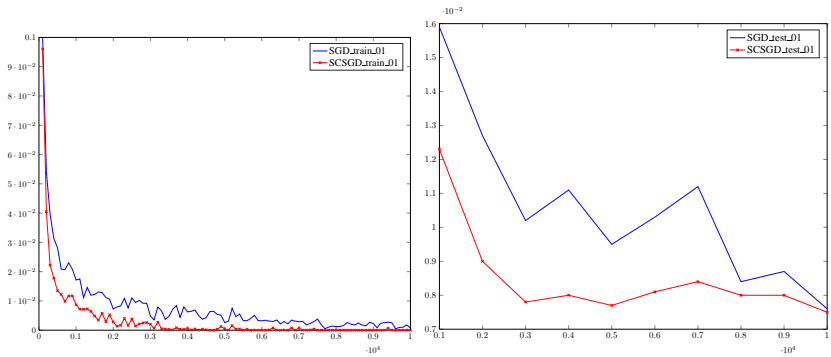

For training, we used the multiclass log loss function.
\figref{fig:mnist} and \figref{fig:mnist01} show the training and the
test errors both w.r.t. the multiclass log loss function and the
zero-one loss (where the $x$-axis corresponds to the number of
iterations). In both graphs, we can see that SCSGD enjoys a much
faster convergence rate.

\begin{figure}
\begin{subfigure} [b] {0.45 \textwidth}
\centering
\resizebox{\linewidth}{!}{% This file was created by matlab2tikz.
% Minimal pgfplots version: 1.3
%  http://www.mathworks.com/matlabcentral/fileexchange/22022-matlab2tikz
\begin{tikzpicture}

\begin{axis}[%
width=6.027778in,
height=4.754167in,
at={(1.011111in,0.641667in)},
scale only axis,
xmin=0,
xmax=10000,
ymode=log,
ymin=0.339111,
ymax=2.24449,
yminorticks=true,
legend style={legend cell align=left,align=left,draw=white!15!black}
]
\addplot [color=blue,solid]
  table[row sep=crcr]{%
100	2.24449\\
200	2.24397\\
300	2.23347\\
400	2.22906\\
500	2.22671\\
600	2.09059\\
700	2.00022\\
800	1.86303\\
900	1.68342\\
1000	1.61515\\
1100	1.48093\\
1200	1.39037\\
1300	1.2763\\
1400	1.17593\\
1500	1.16428\\
1600	1.09533\\
1700	1.09653\\
1800	1.02333\\
1900	0.924695\\
2000	0.891024\\
2100	0.880803\\
2200	0.862393\\
2300	0.852764\\
2400	0.865018\\
2500	0.802697\\
2600	0.821017\\
2700	0.755735\\
2800	0.785847\\
2900	0.726304\\
3000	0.776715\\
3100	0.743903\\
3200	0.708248\\
3300	0.68742\\
3400	0.748628\\
3500	0.708019\\
3600	0.699773\\
3700	0.665745\\
3800	0.693548\\
3900	0.682931\\
4000	0.65608\\
4100	0.597793\\
4200	0.639579\\
4300	0.621286\\
4400	0.640647\\
4500	0.632109\\
4600	0.647567\\
4700	0.607225\\
4800	0.635359\\
4900	0.616159\\
5000	0.622643\\
5100	0.587686\\
5200	0.582771\\
5300	0.617368\\
5400	0.585591\\
5500	0.569579\\
5600	0.57136\\
5700	0.617264\\
5800	0.603717\\
5900	0.594257\\
6000	0.565465\\
6100	0.564583\\
6200	0.55867\\
6300	0.573874\\
6400	0.53808\\
6500	0.540293\\
6600	0.55479\\
6700	0.557227\\
6800	0.545046\\
6900	0.514395\\
7000	0.518723\\
7100	0.545717\\
7200	0.532883\\
7300	0.553391\\
7400	0.535042\\
7500	0.515707\\
7600	0.507584\\
7700	0.554119\\
7800	0.595163\\
7900	0.500718\\
8000	0.516387\\
8100	0.52511\\
8200	0.519758\\
8300	0.530479\\
8400	0.531814\\
8500	0.543282\\
8600	0.527703\\
8700	0.485951\\
8800	0.4747\\
8900	0.48962\\
9000	0.497475\\
9100	0.482804\\
9200	0.473627\\
9300	0.521357\\
9400	0.508806\\
9500	0.568669\\
9600	0.48919\\
9700	0.474951\\
9800	0.434812\\
9900	0.50938\\
10000	0.461923\\
};
\addlegendentry{SGD\_Train};

\addplot [color=red,solid,mark=x,mark options={solid}]
  table[row sep=crcr]{%
100	2.24099\\
200	2.03806\\
300	1.65715\\
400	1.33185\\
500	1.21224\\
600	1.03428\\
700	0.928758\\
800	0.890156\\
900	0.835063\\
1000	0.828314\\
1100	0.748136\\
1200	0.757934\\
1300	0.700956\\
1400	0.671611\\
1500	0.6977\\
1600	0.683744\\
1700	0.652382\\
1800	0.638142\\
1900	0.592987\\
2000	0.601327\\
2100	0.608734\\
2200	0.569193\\
2300	0.580052\\
2400	0.593528\\
2500	0.546472\\
2600	0.586417\\
2700	0.518323\\
2800	0.557667\\
2900	0.52033\\
3000	0.520614\\
3100	0.518412\\
3200	0.538047\\
3300	0.492429\\
3400	0.535209\\
3500	0.525901\\
3600	0.49773\\
3700	0.469681\\
3800	0.543153\\
3900	0.520206\\
4000	0.491571\\
4100	0.430588\\
4200	0.482185\\
4300	0.439448\\
4400	0.443554\\
4500	0.502246\\
4600	0.494962\\
4700	0.445431\\
4800	0.484555\\
4900	0.460889\\
5000	0.488333\\
5100	0.448955\\
5200	0.455718\\
5300	0.465687\\
5400	0.471932\\
5500	0.419828\\
5600	0.463672\\
5700	0.50016\\
5800	0.48526\\
5900	0.451611\\
6000	0.435138\\
6100	0.434012\\
6200	0.434226\\
6300	0.466447\\
6400	0.430329\\
6500	0.419272\\
6600	0.442095\\
6700	0.438184\\
6800	0.440845\\
6900	0.388823\\
7000	0.417332\\
7100	0.423852\\
7200	0.417967\\
7300	0.453016\\
7400	0.401402\\
7500	0.403809\\
7600	0.42453\\
7700	0.435605\\
7800	0.45049\\
7900	0.409383\\
8000	0.41286\\
8100	0.404435\\
8200	0.409273\\
8300	0.428207\\
8400	0.423062\\
8500	0.439372\\
8600	0.426072\\
8700	0.378854\\
8800	0.383803\\
8900	0.390076\\
9000	0.41049\\
9100	0.388279\\
9200	0.396003\\
9300	0.410615\\
9400	0.407064\\
9500	0.453975\\
9600	0.407223\\
9700	0.366613\\
9800	0.339111\\
9900	0.396147\\
10000	0.39052\\
};
\addlegendentry{SCSGD\_train};

\end{axis}
\end{tikzpicture}%
}
\end{subfigure}
\begin{subfigure} [b] {0.45 \textwidth}
\centering
\resizebox{\linewidth}{!}{% This file was created by matlab2tikz.
% Minimal pgfplots version: 1.3
%  http://www.mathworks.com/matlabcentral/fileexchange/22022-matlab2tikz
\begin{tikzpicture}

\begin{axis}[%
width=6.027778in,
height=4.754167in,
at={(1.011111in,0.641667in)},
scale only axis,
xmin=1000,
xmax=10000,
ymode=log,
ymin=0.443682,
ymax=1.58044,
yminorticks=true,
legend style={legend cell align=left,align=left,draw=white!15!black}
]
\addplot [color=blue,solid]
  table[row sep=crcr]{%
1000	1.58044\\
2000	0.973173\\
3000	0.763987\\
4000	0.700481\\
5000	0.636597\\
6000	0.618812\\
7000	0.58456\\
8000	0.560029\\
9000	0.579247\\
10000	0.560449\\
};
\addlegendentry{SGD\_test};

\addplot [color=red,solid,mark=x,mark options={solid}]
  table[row sep=crcr]{%
1000	0.85203\\
2000	0.647555\\
3000	0.566115\\
4000	0.537106\\
5000	0.500949\\
6000	0.490756\\
7000	0.486252\\
8000	0.460926\\
9000	0.466193\\
10000	0.443682\\
};
\addlegendentry{SCSGD\_test};

\end{axis}
\end{tikzpicture}%
%%% Local Variables:
%%% mode: latex
%%% TeX-master: t
%%% End:
}
\end{subfigure}
\caption{SVHN data. Train (left) and test (right) errors w.r.t. the multiclass log loss of SGD and SCSGD}
\label{fig:SVHN}
\end{figure}
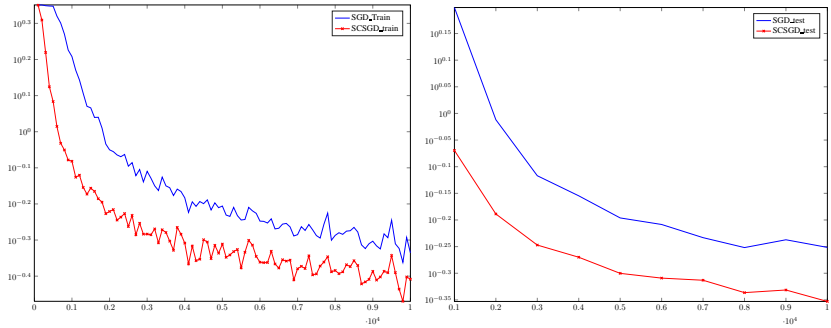

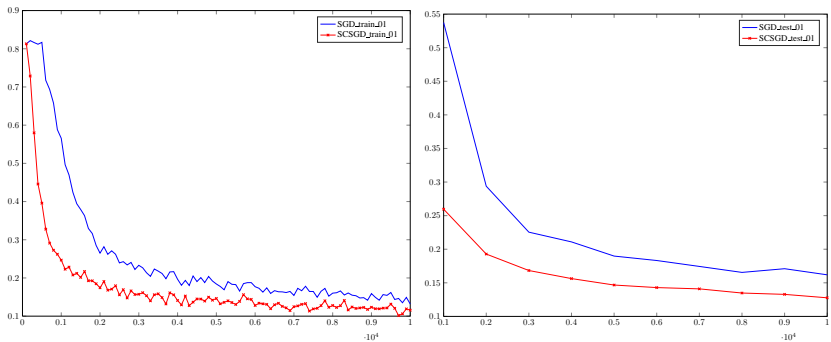
\begin{figure}
\begin{subfigure} [b] {0.45 \textwidth}
\centering
\resizebox{\linewidth}{!}{% This file was created by matlab2tikz.
% Minimal pgfplots version: 1.3
%  http://www.mathworks.com/matlabcentral/fileexchange/22022-matlab2tikz
\begin{tikzpicture}

\begin{axis}[%
width=6.027778in,
height=4.754167in,
at={(1.011111in,0.641667in)},
scale only axis,
xmin=0,
xmax=10000,
ymin=0.1,
ymax=0.9,
legend style={legend cell align=left,align=left,draw=white!15!black}
]
\addplot [color=blue,solid]
  table[row sep=crcr]{%
100	0.812962\\
200	0.821302\\
300	0.816268\\
400	0.81189\\
500	0.816721\\
600	0.717568\\
700	0.693735\\
800	0.657804\\
900	0.588163\\
1000	0.565569\\
1100	0.496487\\
1200	0.469548\\
1300	0.424242\\
1400	0.394277\\
1500	0.379328\\
1600	0.363199\\
1700	0.329535\\
1800	0.316335\\
1900	0.284973\\
2000	0.265027\\
2100	0.281773\\
2200	0.26189\\
2300	0.271016\\
2400	0.262147\\
2500	0.239404\\
2600	0.242662\\
2700	0.234178\\
2800	0.240547\\
2900	0.222227\\
3000	0.233072\\
3100	0.226617\\
3200	0.213519\\
3300	0.204183\\
3400	0.223395\\
3500	0.217947\\
3600	0.211663\\
3700	0.198068\\
3800	0.215703\\
3900	0.216501\\
4000	0.196585\\
4100	0.180798\\
4200	0.193533\\
4300	0.1803\\
4400	0.205113\\
4500	0.190576\\
4600	0.200922\\
4700	0.187927\\
4800	0.203856\\
4900	0.191961\\
5000	0.184192\\
5100	0.177767\\
5200	0.169165\\
5300	0.190403\\
5400	0.183379\\
5500	0.1826\\
5600	0.165625\\
5700	0.185193\\
5800	0.1877\\
5900	0.187787\\
6000	0.177102\\
6100	0.172954\\
6200	0.162863\\
6300	0.174076\\
6400	0.159466\\
6500	0.166716\\
6600	0.163847\\
6700	0.163564\\
6800	0.16183\\
6900	0.164659\\
7000	0.154244\\
7100	0.172756\\
7200	0.166557\\
7300	0.178209\\
7400	0.164647\\
7500	0.164555\\
7600	0.149469\\
7700	0.165093\\
7800	0.172559\\
7900	0.152294\\
8000	0.160238\\
8100	0.161289\\
8200	0.165997\\
8300	0.15541\\
8400	0.160653\\
8500	0.154925\\
8600	0.153427\\
8700	0.147575\\
8800	0.148367\\
8900	0.141617\\
9000	0.159053\\
9100	0.149335\\
9200	0.141698\\
9300	0.156433\\
9400	0.154419\\
9500	0.161836\\
9600	0.143691\\
9700	0.146322\\
9800	0.135132\\
9900	0.149096\\
10000	0.13191\\
};
\addlegendentry{SGD\_train\_01};

\addplot [color=red,solid,mark=x,mark options={solid}]
  table[row sep=crcr]{%
100	0.812962\\
200	0.728712\\
300	0.579857\\
400	0.44591\\
500	0.39587\\
600	0.327794\\
700	0.291241\\
800	0.272506\\
900	0.262025\\
1000	0.246703\\
1100	0.222815\\
1200	0.228451\\
1300	0.207737\\
1400	0.212066\\
1500	0.201469\\
1600	0.216881\\
1700	0.192741\\
1800	0.192641\\
1900	0.184863\\
2000	0.174338\\
2100	0.190703\\
2200	0.167976\\
2300	0.170671\\
2400	0.179587\\
2500	0.155839\\
2600	0.169337\\
2700	0.147572\\
2800	0.166214\\
2900	0.156627\\
3000	0.15722\\
3100	0.161501\\
3200	0.153533\\
3300	0.140109\\
3400	0.156101\\
3500	0.158428\\
3600	0.148712\\
3700	0.131932\\
3800	0.161029\\
3900	0.15586\\
4000	0.141035\\
4100	0.129536\\
4200	0.152668\\
4300	0.127626\\
4400	0.136673\\
4500	0.145069\\
4600	0.144805\\
4700	0.139449\\
4800	0.149681\\
4900	0.141734\\
5000	0.146309\\
5100	0.132399\\
5200	0.136095\\
5300	0.140229\\
5400	0.13568\\
5500	0.13055\\
5600	0.138634\\
5700	0.156231\\
5800	0.145335\\
5900	0.143868\\
6000	0.127907\\
6100	0.133885\\
6200	0.132421\\
6300	0.130815\\
6400	0.119536\\
6500	0.129827\\
6600	0.133693\\
6700	0.125595\\
6800	0.121865\\
6900	0.114513\\
7000	0.124799\\
7100	0.126879\\
7200	0.13088\\
7300	0.13289\\
7400	0.113273\\
7500	0.11909\\
7600	0.120464\\
7700	0.127579\\
7800	0.140192\\
7900	0.123204\\
8000	0.128008\\
8100	0.122716\\
8200	0.127388\\
8300	0.141033\\
8400	0.116231\\
8500	0.123986\\
8600	0.120129\\
8700	0.121553\\
8800	0.122831\\
8900	0.117522\\
9000	0.123238\\
9100	0.119542\\
9200	0.119287\\
9300	0.120968\\
9400	0.121465\\
9500	0.131538\\
9600	0.120529\\
9700	0.101552\\
9800	0.106216\\
9900	0.119063\\
10000	0.115599\\
};
\addlegendentry{SCSGD\_train\_01};

\end{axis}
\end{tikzpicture}%
}
\end{subfigure}
\begin{subfigure} [b] {0.45 \textwidth}
\centering
\resizebox{\linewidth}{!}{% This file was created by matlab2tikz.
% Minimal pgfplots version: 1.3
%  http://www.mathworks.com/matlabcentral/fileexchange/22022-matlab2tikz
\begin{tikzpicture}

\begin{axis}[%
width=6.027778in,
height=4.754167in,
at={(1.011111in,0.641667in)},
scale only axis,
xmin=1000,
xmax=10000,
ymin=0.1,
ymax=0.55,
legend style={legend cell align=left,align=left,draw=white!15!black}
]
\addplot [color=blue,solid]
  table[row sep=crcr]{%
1000	0.5374\\
2000	0.2938\\
3000	0.2254\\
4000	0.2109\\
5000	0.1898\\
6000	0.1832\\
7000	0.1744\\
8000	0.1655\\
9000	0.171\\
10000	0.1619\\
};
\addlegendentry{SGD\_test\_01};

\addplot [color=red,solid,mark=x,mark options={solid}]
  table[row sep=crcr]{%
1000	0.2596\\
2000	0.1928\\
3000	0.1683\\
4000	0.1563\\
5000	0.1467\\
6000	0.1429\\
7000	0.141\\
8000	0.1348\\
9000	0.1328\\
10000	0.1277\\
};
\addlegendentry{SCSGD\_test\_01};

\end{axis}
\end{tikzpicture}%
}
\end{subfigure}
\caption{SVHN data. Train (left) and test (right) errors w.r.t. the zero-one loss of SGD and SCSGD}
\label{fig:SVHN01}
\end{figure}

\paragraph{SVHN:} 
In this experiment we used a much smaller network. The input is images
of size $32 \times 32$ pixels. Using the same terminology as above,
the architecture is now: conv 5x5x8, relu, conv 5x5x16, maxpool 2x2,
conv 5x5x16, maxpool 2x2, affine 32, relu, affine 32, relu, affine 10,
avgpool 4x4.  The results are summarized in the graphs of
\figref{fig:SVHN} and \figref{fig:SVHN01}.  We again see a superior
convergence rate of SCSGD relatively to SGD.

\section*{Acknowledgments} 
This work is supported by the Intel Collaborative Research Institute for Computational Intelligence (ICRI-CI).

\newpage 

\bibliographystyle{plain}
\bibliography{bib}

\begin{thebibliography}{10}

\bibitem{amari1998natural}
Shun-Ichi Amari.
\newblock Natural gradient works efficiently in learning.
\newblock {\em Neural computation}, 10(2):251--276, 1998.

\bibitem{amari2000adaptive}
Shun-Ichi Amari, Hyeyoung Park, and Kenji Fukumizu.
\newblock Adaptive method of realizing natural gradient learning for multilayer
  perceptrons.
\newblock {\em Neural Computation}, 12(6):1399--1409, 2000.

\bibitem{becker1988improving}
Sue Becker and Yann Le~Cun.
\newblock Improving the convergence of back-propagation learning with second
  order methods.
\newblock In {\em Proceedings of the 1988 connectionist models summer school},
  pages 29--37. San Matteo, CA: Morgan Kaufmann, 1988.

\bibitem{bengio2012practical}
Yoshua Bengio.
\newblock Practical recommendations for gradient-based training of deep
  architectures.
\newblock In {\em Neural Networks: Tricks of the Trade}, pages 437--478.
  Springer, 2012.

\bibitem{bordes2009sgd}
Antoine Bordes, L{\'e}on Bottou, and Patrick Gallinari.
\newblock Sgd-qn: Careful quasi-newton stochastic gradient descent.
\newblock {\em The Journal of Machine Learning Research}, 10:1737--1754, 2009.

\bibitem{duchi2011adaptive}
John Duchi, Elad Hazan, and Yoram Singer.
\newblock Adaptive subgradient methods for online learning and stochastic
  optimization.
\newblock {\em The Journal of Machine Learning Research}, 12:2121--2159, 2011.

\bibitem{lecun1998gradient}
Yann LeCun, L{\'e}on Bottou, Yoshua Bengio, and Patrick Haffner.
\newblock Gradient-based learning applied to document recognition.
\newblock {\em Proceedings of the IEEE}, 86(11):2278--2324, 1998.

\bibitem{lecun1998mnist}
Yann LeCun and Corinna Cortes.
\newblock The mnist database of handwritten digits, 1998.

\bibitem{martens2010deep}
James Martens.
\newblock Deep learning via hessian-free optimization.
\newblock In {\em Proceedings of the 27th International Conference on Machine
  Learning (ICML-10)}, pages 735--742, 2010.

\bibitem{moller1993exact}
Martin~F M{\o}ller.
\newblock Exact calculation of the product of the hessian matrix of
  feed-forward network error functions and a vector in 0 (n) time.
\newblock {\em DAIMI Report Series}, 22(432), 1993.

\bibitem{nesterov2004introductory}
Yurii Nesterov.
\newblock {\em Introductory lectures on convex optimization}, volume~87.
\newblock Springer Science \& Business Media, 2004.

\bibitem{netzer2011reading}
Yuval Netzer, Tao Wang, Adam Coates, Alessandro Bissacco, Bo~Wu, and Andrew~Y
  Ng.
\newblock Reading digits in natural images with unsupervised feature learning.
\newblock In {\em NIPS workshop on deep learning and unsupervised feature
  learning}, 2011.

\bibitem{pascanu2013revisiting}
Razvan Pascanu and Yoshua Bengio.
\newblock Revisiting natural gradient for deep networks.
\newblock {\em arXiv preprint arXiv:1301.3584}, 2013.

\bibitem{pearlmutter1994fast}
Barak~A Pearlmutter.
\newblock Fast exact multiplication by the hessian.
\newblock {\em Neural computation}, 6(1):147--160, 1994.

\bibitem{roux2008topmoumoute}
Nicolas~L Roux, Pierre-Antoine Manzagol, and Yoshua Bengio.
\newblock Topmoumoute online natural gradient algorithm.
\newblock In {\em Advances in neural information processing systems}, pages
  849--856, 2008.

\bibitem{saito1997partial}
Kazumi Saito and Ryohei Nakano.
\newblock Partial bfgs update and efficient step-length calculation for
  three-layer neural networks.
\newblock {\em Neural Computation}, 9(1):123--141, 1997.

\bibitem{sarlos2006improved}
Tamas Sarlos.
\newblock Improved approximation algorithms for large matrices via random
  projections.
\newblock In {\em Foundations of Computer Science, 2006. FOCS'06. 47th Annual
  IEEE Symposium on}, pages 143--152. IEEE, 2006.

\bibitem{schraudolph2007stochastic}
Nicol Schraudolph, Jin Yu, and Simon G{\"u}nter.
\newblock A stochastic quasi-newton method for online convex optimization.
\newblock {\em AISTATS}, 2007.

\bibitem{schraudolph2002fast}
Nicol~N Schraudolph.
\newblock Fast curvature matrix-vector products for second-order gradient
  descent.
\newblock {\em Neural computation}, 14(7):1723--1738, 2002.

\bibitem{sutskever2013importance}
I.~Sutskever, J.~Martens, G.~Dahl, and G.~Hinton.
\newblock On the importance of initialization and momentum in deep learning.
\newblock In {\em ICML}, 2013.

\bibitem{vinyals2011krylov}
Oriol Vinyals and Daniel Povey.
\newblock Krylov subspace descent for deep learning.
\newblock {\em arXiv preprint arXiv:1111.4259}, 2011.

\bibitem{werbos1988backpropagation}
Paul~J Werbos.
\newblock Backpropagation: Past and future.
\newblock In {\em Neural Networks, 1988., IEEE International Conference on},
  pages 343--353. IEEE, 1988.

\bibitem{woodruff2014sketching}
David~P Woodruff.
\newblock Sketching as a tool for numerical linear algebra.
\newblock {\em arXiv preprint arXiv:1411.4357}, 2014.

\end{thebibliography}

\newpage
\appendix

\section{Proofs Omitted from The Text}  \label{app:proofs}
\begin{proof} \textbf{(of \lemref{lem:md})}
Using the notation from \secref{sec:intro}, denote 
\[
\Delta_t = \frac{1}{2} D_A(W^\star,W_t)- \frac{1}{2} D_A (W^\star, W_{t+1})  ~.
\]
As in the standard proof of SGD (for the Lipschitz case), we consider the progress of $W_t$ towards $W^\star$. Recall that $W_{t+1} = W_t - \eta \nabla \ell_{y_{i_t}}(W_t x_{i_t})  x_{i_t}^\top A^{-1}$, where $i_t \in [m]$ is the random index that is drawn at time $t$. For simplicity, denote the $p \times n$ matrix $\nabla \ell_{y_{i_t}}(W_t x_{i_t})x_{i_t}^\top$ by $G_t$. Thus, $W_{t+1} - W_t = G_t A^{-1}$. 
By standard algebraic manipulations we have 
\begin{align*}
\Delta_t &=  \frac{1}{2}\tr(A(W^\star-W_{t+1})^\top  (W^\star-W_{t+1})) - \frac{1}{2} \tr(A(W^\star -W_t)^\top (W^\star-W_t)) \\
&= \tr((W^\star-W_{t+1}) A (W_{t+1}-W_t)^\top) + \frac{1}{2}
\tr((W_{t+1}-W_t)^\top A (W_{t+1}-W_t)) \\
&= \tr( (W^\star-W_{t+1}) A ( -\eta G_t A^{-1} )^\top)  + \frac{1}{2}\tr((\eta G_t A^{-1}) A (\eta G_t A^{-1})^\top) \\
& = \eta  \cdot \tr( (W_{t+1} - W^\star)  G_t ^\top) +  \frac{\eta^2}{2}\tr( G_t  A^{-1}  G_t^\top) \\
& = \eta  \cdot \tr( G_t ^\top (W_{t+1} - W^\star) ) +  \frac{\eta^2}{2}\tr( G_t  A^{-1}  G_t^\top) \\
& = \eta  \cdot \tr( G_t ^\top (W_{t} - W^\star) )  + \eta \cdot \tr(G_t^\top (W_{t+1} - W_{t}))+  \frac{\eta^2}{2}\tr( G_t  A^{-1}  G_t^\top) \\
& = \eta  \cdot \tr( G_t ^\top (W_{t} - W^\star) )  - \eta^2 \cdot \tr(G_t A^{-1} G_t^\top)+  \frac{\eta^2}{2}\tr( G_t  A^{-1}  G_t^\top) \\
& = \eta  \cdot \tr( G_t ^\top (W_{t} - W^\star) )  -  \frac{\eta^2}{2}\tr( G_t  A^{-1}  G_t^\top) \\
& \ge  \eta \inner{G_t,W^\star-W_{t+1}} - \frac{\rho^2 \eta^2}{2} \tr(x_{i_t}^\top A^{-1} x_{i_t})~.
\end{align*}
% &= \eta (w_{t+1} - u)^\top \nabla_t   + \frac{1}{2}
% (w_{t+1}-w_t)^\top A (w_{t+1}-w_t)  \\
% &= \eta (w_t - u)^\top \nabla_t - \eta (w_t - w_{t+1})^\top \nabla_t  + \frac{1}{2}
% (w_{t+1}-w_t)^\top A (w_{t+1}-w_t)  \\
% &= \eta (w_t - u)^\top \nabla_t - \eta^2 \nabla_t^\top A^{-1}  \nabla_t  + \frac{1}{2}
% (w_{t+1}-w_t)^\top A (w_{t+1}-w_t)  \\
% &= \eta (w_t - u)^\top \nabla_t - \eta^2 \nabla_t^\top A^{-1}  \nabla_t  + \frac{\eta^2}{2}
% \nabla_t^\top A^{-1} A A^{-1} \nabla_t \\
% &= \eta (w_t - u)^\top \nabla_t - \frac{\eta^2}{2} \nabla_t^\top A^{-1}  \nabla_t  \\
% &\ge \eta (w_t - u)^\top \nabla_t - \frac{\eta^2 L^2}{2} x_t^\top A^{-1}  x_t  ~,
% \end{align*}
where in the last inequality we used the fact that the loss function is
$\rho$-Lipschitz. Summing over $t$ and dividing by $\eta$ we obtain
\[
\sum_{t=1}^T \inner{G_t,W_t-W^\star} ~\le~ \frac{1}{2\eta} \tr(A(W^\star-W_1)^\top
(W^\star-W_1)) + \frac{\eta \rho^2}{2} \tr(A^{-1} \sum_t x_{i_t} x_{i_t}^\top)~.
\]
Recall that $W_1 = 0$. Note that the expected value of $G_t$ is the gradient of $L$ at $W_t$ and the expected value of $x_{i_t} x_{i_t} ^\top$ is $C$. Taking expectation over the choice of $i_t$ for all $t$, dividing by $T$ and relying on the fact that $L(W_t)-L(W^\star) \le \inner{\nabla L(W_t),W_t-W^\star}$, we obtain
\[
L(\bar{W}) - L(W^\star) ~\le~\frac{1}{2 \eta T} \tr(A {W^\star}^\top W^\star) +  \frac{\eta \rho^2}{2} \tr(A^{-1} C)
\]

\end{proof}

 \begin{proof} \textbf{(of \thmref{thm:optBound})}
For simplicity, we assume that $C$ has full rank. If this is not the case, one can
add a tiny amount of noise to the instances to make sure that $C$ is
of full rank. 

We would like to optimize $\tr(A) + \tr(A^{-1}C)$ over all positive definite matrics. Since every matrix $A \succ 0$ can be written as $A = \tau M$, where $M \succ 0$, $\tr(M)=1$ and $\tau = \tr(A)$, an equivalent objective is given by 
\begin{equation} \label{eq:mTrOne}
\min_{\tau>0} \min_{\substack{M \succ 0:\\ \tr(M)=1}} \frac{\sigma^2}{2 \eta T} \tau +\frac{\eta \rho^2}{2 \tau} \tr(M^{-1} C)~.  
\end{equation}
The following lemma characterizes the optimizer. 
\begin{lemma} \label{lem:fanInv}
Let $C \succ 0$. Then,
\[
\min_{\substack{M \succ 0:\\ \tr(M) \le 1}} \tr(M^{-1} C) = (\tr(C^{1/2}))^2~,  
\]
and the minimum is attained by $M^\star=(\tr(C^{1/2}))^{-1} \cdot  C^{1/2}$.
\end{lemma}
Straightforward optimization over $\tau$ yields the value $\tau =
\tr(C^{1/2})$. Subtituitng $\tau$ and $M$ in \eqref{eq:mTrOne} and
applying \lemref{lem:md}, we conclude the proof of \thmref{thm:optBound}. 
\end{proof}
\begin{proof} \textbf{(of \lemref{lem:fanInv})}
First, it can be seen that $M^\star$ is feasible and attains the claimed minimal value. We complete the proof by showing the following inequality for any feasible $A$:
\[
\tr(M^{-1} C) \ge (\tr(C^{1/2}))^2~.
\]
We claim the following analogue of Fan's inequality: For any symmetric matrix $M \in \reals^{n \times n}$, 
\[
\tr(M^{-1} C) \ge \inner{ \lambda^{\uparrow} (M^{-1}), \lambda^{\downarrow} (C)} = \inner{ (\lambda^{\downarrow} (M))^{-1}, \lambda^{\downarrow} (C)} ~,
\]
where $\downarrow$ and $\uparrow$ are used to represent decreasing and increasing orders, respectively, and for a vector $x=(x_1,\ldots,x_n)$ with positive components, $x^{-1}=(1/x_1,\ldots,1/x_n)$. The equality is clear so we proceed by proving the inequality. Let $M$ be a $n \times n$ symmetric matrix. Assume that $C= \sum_{i=1} ^n \lambda_i u_i u_i^\top$ and $M = \sum_{i=1}^n \mu_i v_i v_i^\top$ are the spectral decompositions of $C$ and $M$, respectively. Letting $\alpha_{i,j}=\inner{u_i,v_j}$, we have
\[
\tr(M^{-1} C) = \sum_{i,j} \alpha_{i,j}^2 \lambda_i / \mu_j~.
\]
Note that since both $v_1,\ldots, v_n$ and $u_1,\ldots,u_n$ form orthonormal bases, the matrix $Z \in \reals^{n \times n}$ whose $(i,j)$-th element is $\alpha_{i,j}^2$ is doubly stochastic. So, we have 
\[
\tr(M^{-1} C) = \lambda^\top Z \mu^{-1}~.
\]
Viewing the right side as a function of $Z$, we can apply Birkhoff's theorem and conclude that the minimum is obtained by a permutation matrix. The claimed inequality follows. 

Thus, we next consider the objective 
\[
\min_{\mu \in E} \sum_{i=1} ^n \lambda_i/\mu_i~,
\]
where $E = \{\mu \in \reals^n_+: \sum_{i=1} ^n \mu_i \le 1\}$. The corresponding Lagrangian\footnote{The strict inequalities $\mu_i>0$ are not allowed, but we can replace them with weak inequalities and let $f(\mu)=\infty$ for any $\mu$ whose one of its components is not greater than zero} is  
\begin{align*}
L(\mu; \alpha) = \sum_{i=1} ^n \lambda_i/\mu_i - \sum_{i=1} ^n \alpha_i \mu_i+\alpha_{n+1}(\sum_{i=1} ^n \mu_i-1)~.
\end{align*}
Next, we compare the differential to zero and rearrange, to obtain
\[
(\lambda_i/\mu_i^2)_{i=1}^n = (\alpha_{n+1}-\alpha_i)_{i=1}^n~.
\]
By complementary slackness, $\alpha_1=\ldots=\alpha_n=0$. Thus,
\[
\mu_i ^2 = c \lambda_i^{1/2}~,
\]
for some $c>0$. The constraint $\sum_{i=1} \mu_i \le 1$ implies that $c = \sum_{i=1} \lambda_i ^{1/2}$. Substituting the minimizer in the objective, we conclude the proof.
\end{proof}

\begin{proof} \textbf{(of \lemref{lem:lowCondInv})}
Since $B=EE^\top$ for some matrix $E$, it follows that $QBQ^\top = QE(QE)^\top$, thus it is positive semidefinite. The matrix $a(I-QQ^\top)$ is clearly positive semidefinite. It remains to show that $A$ is invertible and thus it is positive definite. We have
\begin{align*}
A A^{-1} &= (Q B Q^\top + a (I - Q Q^\top)) ( Q B^{-1} Q^\top +
\frac{1}{a} (I-QQ^\top) ) \\
&= Q B Q^\top  Q B^{-1} Q^\top + (I - Q Q^\top) (I - Q Q^\top) + 0 + 0 \\
&= QQ^\top + I - QQ^\top \\
&= I ~.
\end{align*}
 \end{proof}

\begin{proof}  \textbf{(of \thmref{thm:generalLowCond})}
Recall that 
\[
{A} = QBQ^\top + a (I-QQ^\top)~.
\]
According to \lemref{lem:md}, we need to show that 
\[
\tr({A})+\tr({A}^{-1}C) \le \tr(B) + \tr(B^{-1} \tilde{C}) + 2\sqrt{(n-k)(\tr(C) - \tr(\tilde{C}))}~.
\]
Since the trace is invariant to cyclic permutations, we have
\begin{align*}
\tr({A}) &= \tr(QB Q^\top)  + a \cdot \tr(I-QQ^\top) \\&
= \tr(Q^\top Q B)+a(n-k) \\&
= \tr(B) +a(n-k) ~.
\end{align*}
Using \lemref{lem:lowCondInv}, we obtain
\begin{align*}
\tr({A}^{-1}C) &= \tr(Q B^{-1} Q^\top C) + a^{-1} \cdot \tr((I-QQ^\top)C) \\&
= \tr(B^{-1} Q^\top CQ) + a^{-1}(\tr(C) - \tr(QQ^\top C)) \\&
= \tr(B^{-1} \tilde{C}) + a^{-1}(\tr(C) - \tr(\tilde{C})) ~.
\end{align*}
Subtituting $a= \sqrt{\frac{\tr(C)-\tr(\tilde{C})}{n-k}}$, we complete the proof.
\end{proof}

\begin{proof} \textbf{(of \thmref{thm:lowCondExact})}
Note that
\[
\tilde{C} = Q^\top C Q = U_k^\top U D U^\top U_k = U_k D_k U_k) = C_k~.
\]
\[
B = \tilde{C}^{1/2} = C_k^{1/2}~.
\]
\[
B^{-1} \tilde{C} = \tilde{C}^{-1/2} \tilde{C} = \tilde{C}^{1/2} = C_k^{1/2}~.
\]
Invoking \thmref{thm:generalLowCond}, we obtain the desired bound.
\end{proof}

\begin{proof} \textbf{(of \lemref{lem:tildeC_C})}
With a alight abuse of notation, we consider the decomposition $C = C_k + C_{n-k}$ (here, $C_{n-k}$ corresponds to the last $n-k$ eigenvalues rather than to the first $n-k$ eigenvalues).
We need to show that
\[
\tr(\tilde{C}) \ge \tr(C_k) - 3 \tr(C_{d-k})~.
\]
Let $\bar{X} = \frac{1}{\sqrt{m}} X$, where $X \in \reals^{n \times m}$ is the matrix whose columns are $x_1,\ldots, x_m$. Note that $C=\bar{X} \bar{X}^\top$. Also, since  $Q$ satisfies \eqref{eq:woodruff} w.r.t. $X$, it satisfies the same inequality w.r.t. $\bar{X}$. Let $\bar{X}=U\Sigma V^\top$ be the SVD of $X$. Note that the same matrix $U$ participates in the SVD (EVD) of the matrix $C = \bar{X}\bar{X}^\top$, i.e., we have $C = UDU^\top$, where $D=\Sigma^2$. Recall that the best rank-$k$ approximation of $\bar{X}$ is $U_k U_k^\top \bar{X} = U_k \Sigma_k V_k^\top$. By assumption, 
\begin{equation} \label{eq:Q-barX}
\|QQ^\top \bar{X} - \bar{X}\|_F^2 \le 4 \|U_kU_k^\top \bar{X} - \bar{X}\|_F^2
\end{equation}
Note that 
\begin{align*}
\|\bar{X} - QQ^\top \bar{X} \|_F^2 &= \tr(\bar{X}^\top \bar{X}) + \tr(\bar{X}^\top QQ^\top QQ^\top \bar{X}) - 2\tr(\bar{X}^\top QQ^\top \bar{X})\\
&= \tr(C) -\tr(QQ^\top C) = \tr(C) - \tr(\tilde{C})~.
\end{align*}
Similarly, 
\[
\|U_kU_k^\top \bar{X} - \bar{X}\|_F^2 = \tr(C) - \tr(U_kU_k^\top C) = \tr(C) - \tr(C_k) ~.
\]
Thus, \eqref{eq:Q-barX} implies that
\[
\tr(C) - \tr(\tilde{C}) \le 4(\tr(C)- \tr(C_k))~.
\]
Hence,
\[
\tr(\tilde{C}) \ge 4 \tr(C_k) - 3 \tr(C) = \tr(C_k) - 3\tr(C_{d-k})~.
\]
\end{proof}

\begin{proof} \textbf{(of \lemref{lem:everyQ})}
First, we note that $B = B^{-1} \tilde{C} = \tilde{C}^{1/2}$. Thus, we need to show that $\tr(\tilde{C}^{1/2} ) = O( \tr(C_k^{1/2}))$. Second, we observe that for every positive scalar $b$, we have
\[
\tr(\tilde{C}^{1/2} ) = O( \tr(C_k^{1/2}))  \Leftrightarrow \tr(b\tilde{C}^{1/2} ) = O( \tr(bC_k^{1/2}))~.
\]
Denote the $k$ top eigenvalues of $C$ and $\tilde{C}$ by $\lambda_1,\ldots,\lambda_k$ and $\tilde{\lambda}_1,\ldots,\tilde{\lambda}_k$, respectively. According to the above observation, we may assume w.l.o.g. that $\lambda_i \ge 1$ for all $i \in [k]$ (simply consider $b=\lambda_k^{-1}$).

Let $\bar{X} = U \Sigma V^\top$ be the SVD of $\bar{X}$, where $\bar{X} = (1/\sqrt{m})X$. Since $U_k U_k ^\top \bar{X}$ is the best rank-$k$ approximation of $\bar{X}$, we have
\[
\|U_kU_k^\top \bar{X} - \bar{X}\|_F^2 \le \|QQ^\top \bar{X} - \bar{X} \|_F^2
\]
for all $Q \in \reals^{n \times k}$ with orthonormal columns. As in the proof of \lemref{lem:tildeC_C}, this implies that
\[
\tr(C_k) \ge \tr(\tilde{C})~.
\]
Therefore,
\begin{align*}
\tr(\tilde{C}^{1/2}) - \tr(C_k^{1/2}) &= \sum_{i=1} ^k (\sqrt{\tilde{\lambda}_i} - \sqrt{\lambda_i} )
=\sum_{i=1} ^k \frac{\tilde{\lambda}_i - \lambda_i}{\sqrt{\tilde{\lambda}_i} + \sqrt{\lambda_i}} \\
&\le \sum_{i=1} ^k \tilde{\lambda}_i - \lambda_i \\
&=\tr(\tilde{C}) - \tr(C_k)) \le 0
\end{align*}
where the first inequality follows from the assumption that $\lambda_i \ge 1$ for all $i\ \in [k]$.
\end{proof}

\end{document}